\documentclass[12pt]{article}

\usepackage{amsmath}
\usepackage{amssymb}
\usepackage{enumerate}
\usepackage{theorem}
\usepackage{amsfonts}
\usepackage[a4paper,margin=2cm,vmargin={1cm,2cm},includeheadfoot]{geometry}
\usepackage{ifthen}
\usepackage[scaled=0.92]{helvet}
\usepackage[latin1]{inputenc}
\usepackage{srcltx}
\usepackage{tikz}
\usepackage[labelfont=bf,margin=1cm,justification=justified,labelsep=period]{caption}

\definecolor{gris}{gray}{0.8}
\newcommand{\R}{{\mathbb {R}}}
\newcommand{\N}{{\mathbb {N}}}
\newcommand{\Q}{{\mathbb{Q}}}
\newcommand{\D}{{\mathbb{D}}}
\newcommand{\1}{{\bf 1}}

\newcommand{\nin}{\noindent}
\newcommand{\wt}{\widetilde}
\newcommand{\wh}{\widehat}
\newcommand{\buil}[3]{\mathrel{\mathop{\kern0pt#1}\limits_{#2}^{#3}}}
\newcommand{\prob}{\mathbb{P}}
\newcommand{\bs}{\boldsymbol}
\newcommand{\eps}{{\varepsilon}}        % epsilon

\theoremheaderfont{\bfseries}
\newtheorem{theorem}{Theorem}
\newtheorem{definition}[theorem]{Definition}
\newtheorem{proposition}[theorem]{Proposition}
\newtheorem{lemma}[theorem]{Lemma}

\newtheorem{corollary}[theorem]{Corollary}
{\theorembodyfont{\rmfamily}
\newtheorem{remark}[theorem]{Remark}
\newtheorem{remarks}[theorem]{Remarks}}

\newcommand{\E}{{\mathbb{E}}}

\newcommand{\Prob}{\mathbb{P}}
\numberwithin{equation}{section} \numberwithin{theorem}{section}
\begin{document}

\title{Local Malliavin Calculus for L\'{e}vy Processes and Applications
\thanks{Partially supported by the  Ministerio de Ciencia e
Inovaci\'on grant MTM2009-08869 and FEDER}}
\date{08-03-2012}

\author{Jorge A. Le\'{o}n\thanks{Partially supported by the CONACyT grant
98998}\\
Departamento de Control Autom\'{a}tico\\
Cinvestav-IPN\\
Apartado Postal 14-740\\ 07000 M\'{e}xico D.F.,
Mexico\\jleon@ctrl.cinvestav.mx \and
 Josep L. Sol\'{e}, Frederic Utzet
\\ Departament de
Matem\`{a}tiques\\ Universitat Aut\`{o}noma de Barcelona\\
08193-Cerdanyola (Catalunya), Spain\\jllsole@mat.uab.cat\\
utzet@mat.uab.cat \and Josep Vives\\ Departament de Probabilitat,
L\`{o}gica i Estad\'{\i}stica\\ Universitat de Barcelona\\
08007-Barcelona (Catalunya), Spain\\josep.vives@ub.edu}

\maketitle

\begin{abstract}
%\textbf{Abstract}

%\begin{quote}
 {\small
In this paper  a  Malliavin calculus for L\'{e}vy processes based on
a family of true derivative operators is developed. The starting
point is an extension to Lévy processes of  the pioneering paper by
Carlen and Pardoux \cite{Ca-Pa-90} for the Poisson process, and
 our
approach includes also the classical Malliavin derivative for
Gaussian processes. We obtain a sufficient condition for the
absolute continuity of functionals of the Lévy process. As an
application, we analyze the absolute continuity of the law of the
solution of some stochastic differential equations.}
\end{abstract}
%\medskip

%\end{quote}

\section{Introduction}

  Over the past twenty five years, it has been done several attempts  for
extending the Malliavin calculus to a scenario driven by processes
with jumps. Among these efforts, we should  cite  the celebrated
papers of Bismut \cite{Bis}, L\'eandre \cite{Lean} and the monograph
of   Bichteler, Gravereaux and Jacod \cite{BGJ}. In the latter it is
assumed  the restriction that the L\'evy measure $\nu$ is absolutely
continuous with respect to the Lebesgue measure and the main idea is
to derive with respect to the size of the jumps. However,   if the
measure $\nu$ is discrete, the condition in \cite{BGJ}  is not
satisfied, and by this reason Carlen and Pardoux \cite{Ca-Pa-90}
proposed, for  a Poisson process, a derivative operator based on
translation of the jump times.

There is another approach to define a Malliavin derivative as an
annihilation operator    through a chaos expansion. In this way,
Nualart and Vives \cite{NV} proved that for the standard Poisson
process  the annihilation operator (Malliavin derivative) agrees
with a translation operator. Immediately it appeared the necessity
of extending these methods to L\'evy processes with a more general
L\'evy measure, and many authors proposed and  studied such
difference operator: see for example, Picard \cite{Pic1} and
\cite{Pic2}, members of the norwegian school of probabilities like
Di Nunno, Meyer-Brandis, {\O}ksendal and Proske  \cite{Nunno-04} and
Benth and Lokka \cite{BN}, and also Ishikawa and Kunita \cite{IK},
and   Sol\'e, Utzet and Vives \cite{SUV}.  It is also interesting to
cite the extension to abstract Fock spaces by Applebaum
\cite{Apple}. Note that in all those works, in contrast with the
Wiener case,  the annihilation operator is not  a true derivative
and  it does not satisfy the chain rule.

 A number of authors
have made contributions in the direction of defining  a true
derivation operator for jump processes, under more or less
restrictive assumptions. We should cite   Denis \cite{Denis00} and
Bouleau and Denis \cite{Bou-Denis}, who use the context of Dirichlet
forms, Di Nunno \cite{Nunno-02}, and recently Kulik \cite{kulik08}
and Bally \cite{Bal1}. We have to note  that Decreusefond and Savy
\cite{De-Sa-06} worked  with marked point processes (i.e., processes
with compensator $\lambda(s)ds\nu(dz)),$ introducing an operator
similar to the one studied in this paper, but with the restriction
that $\nu$ should be a probability measure.

In this paper we define a true derivative operator for functionals
of a general L\'evy processes which can have a continuous part. Such
operator is an extension to L\'evy processes of the one introduced
by Carlen and Pardoux \cite{Ca-Pa-90} for the Poisson process. Our
main finding is that this operator can be restricted to a Borel
subset of $\R$, say $\Lambda$, where  live part of the jump sizes of
the Lévy processes, and a suitable (time and space) weight function
$k(t,x)$ can be introduced. That subset and weight function can be
chosen in every application to suit the properties of the problem in
hand. By this reason we call it the {\it local Malliavin
derivative.} Our procedure is, as in the Brownian Malliavin
calculus, to start with a class of smooth random variables, define
the derivative, and extend by density to a domain contained in
$L^2(\Omega)$. We observe that this approach can be used even for
the difference operator as it is done in  Geiss and Laukkarinen
\cite{GL}, but, as we commented, unfortunately, the chain rule is
not satisfied for the difference operator.

We would like to point out that we don't present here a complete
Malliavin calculus theory, because our integration by parts formula
(Proposition \ref{pro:duality}) does not cover the general case. On
the other hand, we do not give any result about the smoothness of
the density. However we believe that we may obtain some interesting
results from our development.

Using the chain rule of  our  derivative operator     we give a
sufficient condition for the absolutely continuity of functionals of
the involved L\'evy process analogue  to the one of Carlen and
Pardoux \cite{Ca-Pa-90}. In particular,  we  apply that condition to
deduce the absolute continuity of  solutions of different stochastic
differential equations. The first case studied is the equation
governed  by a L\'evy process with continuous part
\begin{equation}\label{eq:1.1}
Z_t=x+\int_0^t
b(Z_s)ds+\int_0^t\sigma(Z_s)dW_s+\int_0^t\int_{\R_0}l(y)ydN(s,y),\quad
t\in[0,T],
\end{equation}
 with  suitable coefficients $b,\sigma$ and $l,$ and  where
$\mathbb{R}_0:={\mathbb{R}}-\{0\}$. We get a similar result as the
one contained in   Nualart \cite{Nu-95}, where jumps are not
considered.

The second application is devoted to the following equation, driven
by a pure additive jump process,
\begin{equation}\label{eq:1.2}
Z_t=x+\int_0^tf(Z_s)ds+\int_0^t\int_{\R_0}h(y)dN(s,y),\quad
t\in[0,T].
\end{equation}
Here $x\in\R$, $f:\R\rightarrow\R$ is a function with a bounded
derivative and $h\in L^2(\R_0,\nu)\cap L^1(\R_0,\nu)$,   $\nu$ is
the Lévy measure corresponding to the driving Lévy process. It is
well-known that this equation has a unique square-integrable
solution. The first step  is to  find the derivative of the solution
at time $T$.  Then, in the case that $\nu(\R_0)<\infty$  and  that
the function $f$   is monotone, choosing the convenient $k(t,x)$, we
obtain that $Z_T$ satisfies the criterion  established  for the
absolute continuity
 in the set of paths where there is at least one jump.
Moreover, if $\nu(\R_0)=\infty,$ it is only necessary the condition
of monotony for $f$ in a neighborhood of the initial condition $x$.
This is the same type of result  obtained by Nourdin and Simon
\cite{Nu-Si-06} using the
  stratification method due to Davydov {\it et al.} \cite{DLS}. Note that as the L\'evy measure is not  finite, we have, with probability
  one,
  infinity jumps in any neighborhood of the origin, and this fact
guarantees that we can find the suitable function $k$  in order to
restrict us to the neighborhood of monotony of $f$.

 Finally, we  study the equation
\begin{equation}\label{eq:1.3}
Z_t=x+\int_0^tf(Z_s)ds+\int_0^t\int_{\R_0}h(y)g(Z_{s-})dN(s,y),\quad
t\in[0,T],
\end{equation}
with some conditions on $f, \, g$ and $h$. The existence of a
density for $Z_T$ was analyzed  by Carlen and Pardoux
\cite{Ca-Pa-90} in the case that the involved L\'evy process is a
Poisson process. When the Lévy process is a compound Poisson, using
our criterion of absolute continuity with
 a suitable function $k(t,x)$ that allows us to concentrate in
only one jump, we
   prove the absolutely continuity of the solution in the
 set of paths where
 there is at least one jump in the case where
   a condition on the Wronskian of $f$ and $g$ is satisfied.

 The paper is organized as follows. Section 2 presents the
 preliminaries and notations. Section 3
 is devoted to the study of the family of  derivative
  operators introduced and to  the criterion of absolute continuity. In  Section 4 we present  some properties of the derivative operator
  applied  to  multiple integrals which are needed in the next section.
   Finally, in Section 5 we
     get  the absolute
 continuity of the solutions  to the stochastic differential
 equations (\ref{eq:1.1}), (\ref{eq:1.2}) and (\ref{eq:1.3}).

\section{Preliminaries}

In this section we introduce the framework and the notation that we
use in this paper. Let $X=\{X_t, t\geq 0\}$ be a L\'{e}vy process,
that is, $X$ has stationary and independent increments, is
continuous in probability, is c\`{a}dl\`{a}g  and $X_0=0$, with
triplet $(\gamma, \sigma, \nu),$ and it is  defined on a complete
probability space $(\Omega, \mathcal{F}, {\mathbb{P}})$; for all
these notions we refer to Sato \cite{Sa-99}. By convention, the Lévy
measure $\nu$ is defined on $\R$ and $\nu(\{0\})=0$. Also, let
$\{\mathcal{F}_t, t\geq 0\}$ be the natural filtration associated
with $X$. The process $X$ admits the L\'{e}vy-It\^{o}
representation:

\begin{equation}\label{eq:lev-pr}
X_t=\gamma t+\sigma W_t+\int\int_{(0,t]\times\{|x|>1\}} x
dN(s,x)+\lim_{\epsilon \downarrow 0}
\int\int_{(0,t]\times\{\epsilon<|x|\leq 1\})} x d{\tilde N}(s,x),
\end{equation}
where $W$ is a standard Brownian motion, $N$ is the jump measure of
the process defined below, $W$ and $N$ are independent,  and $d{\tilde
N}(t,x)=dN(t,x)-dt \nu (dx).$ The limit is a.s., uniform in $t$ in
any bounded interval.

For any Borel set $B\in\mathcal{B}((0,\infty)\times \mathbb{R}_0)$,
with $\mathbb{R}_0:={\mathbb{R}}-\{0\},$ the jump measure is defined
by
\begin{equation*}
N(B)=\mbox{card}\{t\in (0,\infty): (t, \Delta X_t)\in B\},
\end{equation*}
where $\Delta X_t=X_t-X_{t-}$.

It is well known (see Itô  \cite{It-56}) that  $X$ determines
 a centered  independent scattered  random measure $M$  on
$[0,\infty)\times \mathbb{R}$. To this end,   consider the measure
$$\mu(dt,dx)=\lambda(dt)\delta_0(d x)+\lambda(dt) x^2\nu(dx),$$
where $\lambda$ is the Lebesgue measure in $\R$. For
$E\in\mathcal{B}([0,\infty)\times \mathbb{R})$ such that
$\mu(E)<\infty,$
\begin{equation*}
M(E)=\sigma \int_{E(0)} dW_t+ \lim_{n\rightarrow\infty} \int\int_{\{(t,x)\in E^{\prime}: \frac{1}{n}%
\leq|x|\leq n\}} x d{\tilde N}(t,x),
\end{equation*}
where $E(0)=\{t>0: (t,0)\in E\}$, $E^{\prime}=E-\{(t,0)\in E\}$, and
the limit is in $L^2(\Omega)$. For $A, B\in
\mathcal{B}([0,\infty)\times \mathbb{R})\big)$, with $\mu(A)<\infty$
and $\mu(B)<\infty$, we have
$$\E[M(A)M(B)]=\mu(A\cap B).$$

Because we  work in a finite time interval $[0,T]$,  we consider the
restriction of $M$ to $[0,T]\times \R$. For a function $h\in
L^2([0,T]\times {\mathbb R}, \mu),$ we can construct the integral
(in the $L^2(\Omega)$ sense)
$$M(h):=\int_{[0,T]\times \R}h(t,x)\, dM(t,x),$$
which is  the multiple integral of order 1 defined by  Itô
\cite{It-56}. This integral is centered, and for $g,h\in
L^2([0,T]\times {\mathbb R}, \mu),$
\begin{equation*}
E[M(h)M(g)]=\int_{[0,T]\times \mathbb R}h\, g \, d\mu.
\end{equation*}
We can write  this integral as
\begin{equation*}
M(h)=\sigma \int_0^{T} h(t,0) dW_t+ \int_ {[0,T]\times \R_0}h(t,x)
x\, d{\tilde N}(t,x).
\end{equation*}

In this paper we also need to consider integrals with respect to the
 non compensated counting measure, defined by

\begin{equation}
\label{non-comp}\int_ {[0,T]\times \R_0}h(t,x) x\, d{ N}(t,x):=\int_
{[0,T]\times \R_0}h(t,x) x\, d{\tilde
N}(t,x)+\int_0^T\int_{\R_0}h(t,x)x\, dtd\nu(x),
\end{equation}
when both integrals in the right hand side make sense.

Observe that if $X$ has continuous trajectories, then $M$ is nothing
more than the independent Wiener measure, and if $\sigma=0$ and
$\nu=\delta_{1}$ we obtain the  standard  random Poisson measure.

\section{The derivative operators}\label{sec:3}
\subsection{Carlen-Pardoux derivative}

The starting point of this paper is the following remark due to
Le\'{o}n and Tudor \cite{Le-Tu-98} about the Carlen--Pardoux
derivative \cite{Ca-Pa-90}:
 Let $N$ be a Poisson process of parameter 1, and $\wt N_t=N_t-t$ be the compensated
Poisson process. Denote by $T_1<\cdots <T_n\dots$ the jump times of
$N$. Let $h\in{\cal C}^1\big([0,T]\big).$
  Then, the random variable
  $$Y=\int_0^T h(s)\, d\wt N_s=\sum_{T_n\le T} h(T_n)-\int_0^T h(s)\, ds$$
  is in the domain of the Malliavin derivative in the Carlen-Pardoux sense and the derivative satisfies the following equality:
  $$D_tY=\int_0^T h'(s)\big(\frac{s}{T}-1_{(t,T]}(s)\big)\, dN_s.$$
  In the next subsection, we extend formally this property to define a Malliavin derivative for a Lévy process.

  \subsection{The local Malliavin derivative}

Henceforth,  $\mathcal{S}$ denotes the family of all  functionals of
the form
\begin{equation}\label{smoothfunctional}
f(M(h_1),\dots, M(h_n)),
\end{equation}
where $f$ is in $C^{\infty}_b ({\mathbb R}^n)$ ($f$ and all its
partial derivatives are bounded), $h_1,\dots, h_n\in
L^2([0,T]\times{\mathbb R}, \mu)$ and, for all $x\neq 0$,
$h_i(\cdot,x)\in C^1 ([0,T])$,
 $i\in \{1,\ldots,n\}$,
% such that $h_i$
% and the derivatives have a continuous extension to $[0,T]$,
 and   $\partial h_i\in L^2([0,T]\times{\mathbb R}, \mu)$
where  $\partial$ means the partial derivative with respect to time.
The set $\cal S$ is called the family of smooth random variables,
and we will prove later that it is dense in $L^2(\Omega)$ (see
Proposition \ref{densitat}).

We will also consider the family $\cal K$ of all  bounded functions
$k:[0,T]\times\R_0\to \R$ such that they and their partial
derivative with respect time are in
$L^2([0,T]\times\R_0,\lambda\times\nu)\cap
L^1([0,T]\times\R_0,\lambda\times\nu)$.

Now we introduce the main tool of this paper. Namely, the family of
derivative operators for functionals of the L\'evy process $X$.

\begin{definition}\label{localderivative}
Given  $k\in {\cal K}$, $F\in {\mathcal S}$  and $\Lambda\in
\mathcal{B}(\R)$, we define
\begin{equation*}
D^{\Lambda, k}_{t} F=\sum_{i=1}^n (\partial_i f)(M(h_1),\dots,
M(h_n)) D^{\Lambda, k}_{t}M(h_i), \quad t\in[0,T],
\end{equation*}
where
\begin{equation}\label{eq:dermh}
D^{\Lambda, k}_{t}M(h)=\1_{\Lambda} (0)  \,\sigma
h(t,0)+\int_{[0,T]\times (\Lambda\cap \R_0)} k(s,y)\partial_s h(s,y)
(\frac{s}{T}-\1_{[t,T]}(s))y\, dN(s,y).
\end{equation}
\end{definition}
We call $D^{\Lambda, k}_{t} F$ the local Malliavin derivative of
$F$.

\begin{remarks}\label{re:der}

\begin{enumerate}
\item The integral in the right-hand side of
(\ref{eq:dermh}) is well-defined. Indeed, in agreement with
(\ref{non-comp}),
$$
 \int_{[0,T]\times (\Lambda\cap \R_0)}
k(s,y)\partial_s h(s,y) (\frac{s}{T}-\1_{[t,T]}(s)) y\, d{\tilde
N}(s,y)$$ is finite because $k$ is bounded, and
$$\int_0^T \int_{\Lambda\cap \R_0} k(s,y)\partial_s
h(s,y) (\frac{s}{T}-\1_{[t,T]}(s))y \, ds \nu(dy)$$ is finite by
Cauchy-Schwarz inequality. Moreover, $D^{\Lambda, k}_{t} F\in
L^2([0,T]\times\Omega, \lambda\otimes\Prob)$.
 This explain
why we need to introduce the family ${\cal K}$ in the definition of
the derivative operator. Also in Section \ref{sec:4}, we  see that
this family ${\cal K}$ allows to consider applications to stochastic
differential equations.

\item The function $k(t,x)$ is an  essential ingredient of the derivative, and it may change from one application to another.
To shorten the notation, in general  we will omit that $k$ in the
expression $D^{\Lambda, k}_{t} F$.

\item If $h, \partial h\in L^2([0,T]\times\R_0,\lambda\otimes\nu)$, then
$h(s,y)/y$ and $\partial_s h(s,y)/y$ belong to
$L^2([0,T]\times\R_0,\mu)$. Hence,
$$D^{\Lambda}_{t}M(h(s,y)/y)=\int_{[0,T]\times(\Lambda\cap \R_0)} k(s,y)\partial_s h(s,y)
(\frac{s}{T}-\1_{[t,T]}(s))dN(s,y),$$ where the factor $y$ has
disappeared in the integral.
\item
Observe that if $\Lambda=\{0\}$, then $D^{\Lambda}$ is the Malliavin
derivative with respect to the Brownian part of the L\'evy process
$X$.
\item
If $\Lambda=\{x\}$ for some $x\neq 0$ with $\nu(\{x\})\neq 0$, we
obtain
\begin{equation*}
D^{\{x\}}_{t} M(h)=x  \int_0^T  k(s,x)\partial_s h(s,x) (\frac{s}{T}%
-\1_{[t,T]}(s))dN_s^x,
\end{equation*}
where $N^x_s$ is the Poisson process on $[0,T]$ that counts the
number of jumps of height $x$. Moreover, if the L\'{e}vy process is
the standard Poisson process,
 $x=1$, and if we take $k(t,x)$  independent of the time parameter,   we
obtain
\begin{equation*}
D^{\{1\},k}_{t} M(h)= k(1)\int_0^T \partial_s h(s,1)
(\frac{s}{T}-\1_{[t,T]}(s))dN_s^1.
\end{equation*}
Like we commented in the previous subsection, that  is exactly a
constant times the derivative defined by  Carlen and  Pardoux in
\cite{Ca-Pa-90}, as it is shown in León and Tudor \cite{Le-Tu-98}.

\item The integration by parts formula allow us to write
\begin{eqnarray*}\lefteqn{
 \int_{\Lambda\cap\R_0} \int_0^T k(s,y)
 \partial_s h(s,y) y (\frac{s}{T}-\1_{[t,T]} (s))ds \nu(dy)
}\\
&\quad\quad\phantom{holla}&=\int_{\Lambda\cap\R_0}\big\{k(t,y)h(t,y)
-\frac{1}{T}\int_0^Tk(s,y)
h(s,y)ds\big\}y\nu(dy)\\
&&\ \ \ -\int_0^T\int_{\Lambda\cap\R_0}
h(s,y)\partial_sk(s,y)\left(\frac{s}{T}-\1_{[t,T]}(s)\right)y\nu(dy)ds.
\end{eqnarray*}
Thus,
\begin{eqnarray*}
D^{\Lambda}_{t} M(h)&=&\1_{\Lambda} (0)  \,\sigma h(t,0) +
\int_{[0,T]\times\Lambda\cap\R_0}k(s,y)
\partial_s h(s,y) y (\frac{s}{T}-\1_{[t,T]} (s))d{\tilde
N}(s,y)\\
&&+\int_{\Lambda\cap\R_0}\big\{k(t,y)h(t,y)-\frac{1}{T}\int_0^T
k(s,y)h(s,y)ds\big\}y\nu(dy)\\
&&-\int_0^T\int_{\Lambda\cap\R_0}
h(s,y)\partial_sk(s,y)\left(\frac{s}{T}-\1_{[t,T]}(s)\right)y\nu(dy)ds.
\end{eqnarray*}

\item
If $0\notin\Lambda$, then  $D^{\Lambda}M(h)$ is orthogonal to
the constant functions. In fact, the Fubini theorem (see Proposition
\ref{fubini}) yields
\begin{eqnarray*}\lefteqn{
\int_{0}^{T}\int_{[0,T]\times\Lambda }k(s,y)\partial_s h(s,y)y(\frac{s}{T}%
-\1_{[0,s]}(t))dN(s,y)dt}\\
&&\quad
=\int_{[0,T]\times\Lambda }k(s,y)\partial_s h(s,y)y\int_{0}^{T}(\frac{s}{T}%
-\1_{[0,s]}(t))dtdN(s,y)=0.
\end{eqnarray*}

\end{enumerate}
\end{remarks}

The next result is a Fubini type proposition. Its  proof is
straightforward and omitted.
\begin{proposition}\label{fubini} Let $f\in L^2([0,T]^2\times \R, \lambda^{\otimes 2}\otimes \mu)$.
%$f:[0,T]^2\times \R\to \R$ be a  measurable function such that
%$$\int_{[0,T]}\int_{[0,T]\times \R}f^2(t,s,x)\, dt d\mu(s,x)<\infty.$$
Then the following integrals are well defined and

\begin{equation}
\int_0^TM(f(t,\cdot))\, dt=M(\int_0^Tf(t,\cdot)\,dt),\ \text{a.s.}
\end{equation}
\end{proposition}
\bigskip

\bigskip

In the Brownian Malliavin calculus   if $F$ is a smooth random
variable such that $F=0$, it follows  that $DF=0$. However, for
L\'{e}vy processes matters are different. For example, if $X$ is a
standard Poisson process and we take $h=1$, then
$F=f(M(h))=f(N_T-T)$, and $F=0$ only implies that $f$ is zero on the
set $\N-T$. In  the next proposition we show the corresponding
property.

\begin{proposition}
\label{definit}
 Let $F\in {\cal S}$ such that $F=0 \,  \text{a.s.}$
Then $D^{\Lambda}F=0, \ \lambda\otimes \Prob-\text{a.e.}$
\end{proposition}

\medskip

\nin{\it Proof.} Let $F=f(M(h_1),\dots, M(h_n)).$   To consider the
most general case, assume that $0\in\Lambda$. By definition,
\begin{align*}
&D^{\Lambda}_{t} F=\underbrace{\sum_{i
}(\partial_i f)(M(h_1),\dots,
M(h_n))\sigma h_i(t,0)}_{(*)} \\
&+ \underbrace{\sum_{i=1}^n(\partial_i f)(M(h_1),\dots,
M(h_n))\int_{[0,T]\times (\Lambda\cap \R_0)} k(s,y)\partial_s
h_i(s,y) (\frac{s}{T}-\1_{[t,T]}(s))y\, dN(s,y)}_{(**)}.
\end{align*}

The term (*) is a Malliavin derivative with respect to the Brownian
part (see Solé {\it et al.} \cite{SUV}). From the independence
between the Brownian part and the jumps part of $X$, by fixing the
random variables corresponding to the jumps part, we have that
$$M(h_j)=\sigma W(h_j)+C_j,$$
where $W(h)=\int_0^Th(t,0)\, dW_t$ and $C_j$ is a constant; changing
the function $f$ we can take $C_1=\cdots=C_n=0$.
 First assume that the Gaussian vector
$(W(h_1),\dots,W(h_n))$ has a regular covariance matrix. Then it has
density, and
$$f(W(h_1)\dots, W(h_n))=0, \ \text{a.s.}$$
implies
$$f(x_1,\dots,x_n)=0,\ \lambda^n-\text{a.e.}$$
It follows that all derivatives of $f$ are zero.

If the covariance matrix of $(W(h_1),\dots,W(h_n))$ is singular,
denote by $k$ the rank of that matrix, and (reordering) assume that
$(W(h_1),\dots,W(h_k))$ has a density. Then there are numbers
$\mu_i^j$, $i=1,\dots,n$, $j=k+1,\dots, n$ such that
$$W(h_j)=\sum_{i=1}^k\mu^j_i W(h_i),\ j=k+1,\dots, n.$$
Given the linearity of Itô-Wiener integral, it follows that
$\lambda$--a.s.,
$$h_j=\sum_{i=1}^k\mu^j_i h_i,\ j=k+1,\dots, n.$$
Consider the vectorial function
$$\bs \Phi=(\Phi_1,\dots, \Phi_n):\R^k\to\R^n$$
given  by
$$\Phi_j(x_1,\dots,x_k)=x_j,\ \text{for}\ j=1,\dots, k,$$
and
$$\Phi_j(x_1,\dots,x_k)=\sum_{i=1}^k\mu^j_i x_i,\ j=k+1,\dots, n.$$
Define
$g:=f\circ \bs \Phi.$
 Then
$$f(W(h_1)\dots, W(h_n))=g(W(h_1)\dots, W(h_k)).$$
Since  $\big(W(h_1)\dots, W(h_k)\big)$ has a density,  we have that
$$\frac{\partial g}{\partial x_j}=0, \ j=1,\dots, k.$$
By the chain rule,
$$J_g=J_f \, J_{\bs \Phi},$$
where $J_m$ denotes the Jacobian matrix of a generic function $m$. Writing
$\bs h=(h_1,\dots, h_k)^T$, we have
$$0=J_g \bs h(t)= J_f  J_{\bs \Phi} \bs h(t)=\sum_{j=1}^n \frac{\partial f}{\partial x_j}\, h_j(t),$$
and hence (*) is zero.

\medskip

The term (**) is treated as follows.  Again by the independence
between the Brownian and the jumps part, we can assume $\sigma=0$.
For the sake of simplicity we consider the case $n=2$.

Consider the set
$$\Lambda_m=\Lambda\cap\{x\in \R: 1/m<\vert x\vert <m\},$$
 and denote  by $N^{\Lambda_m}$ the restriction of the counting measure of
jumps $N$ to the jumps with size in $\Lambda_m.$ That is, for $B\in
{\cal B}\big((0,\infty)\times \R_0 \big)$,
$$N^{\Lambda_m}(B)=\#\{t:\ (t,\Delta X_t)\in B \ \text{and} \ \Delta X_t\in \Lambda_m\},$$
and
$$N_T^m=N^{\Lambda_m}([0,T]\times \R_0)<\infty,\ \text{a.s.},$$
due to $\nu(\Lambda_m)<\infty$. The jumps of $X$ in $[0,T]$ with
height in $\Lambda_m$ can be ordered, $T_1<T_2<\cdots< T_{N_T^m}$
Given $N_T^m=k$, the vector   $(T_1,\dots,T_k)$ has a  Dirichlet
distribution on $[0,T]$, independent of the height of the jumps. The
jump size has a  distribution given by $P\{\Delta X_{T_j}\in
A\}=\nu(A)/\nu(\Lambda_m)$, $A\in {\cal B}(\Lambda_m)$.

In a similar way, define $N^{\R_0-\Lambda_m}$. Since $N^{\Lambda_m}$
and $N^{\R_0-\Lambda_m}$ are independent, denoting by
$\E_{\Lambda_m}$ and $\E_{\R_0-\Lambda_m}$ the expectation with
respect to each jump counting measure, and as $F=0$ implies
$F^2=f^2(M(h_1),\dots, M(h_n))=0,$ we have that
%$$\E\big[F\1_{\{N_T^m=k\}}\big]=\E\big[F^2\1_{\{N_T^m=k\}}\big]=\E_{\R_0-\Lambda_m}\E_{\Lambda_m} \big[F\1_{\{N_T^m=k\}}\big]=\E_{\R_0-\Lambda_m}\E_{\Lambda_m} \big[F^2\1_{\{N_T^m=k\}}\big]
%=0.$$
$$0=\E\big[F^2\big]=\E_{\R_0-\Lambda_m}\E_{\Lambda_m} \big[F^2\big]
=\E_{\R_0-\Lambda_m}\Big[\sum_{k\ge 0}\prob\{N^m_T=k\}\E\big[F^2\, |\, N_T^m=k\big]\Big].$$
Due to the fact that $F^2\ge 0$, for each $k$,
\begin{align*}
0&=\E_{\R_0-\Lambda_m}\Big[
\int_{S_k(\Lambda_m)}f^2\big(\sum_{j=1}^kh_1(s_j,y_j)y_j-\int_0^T\int_{\Lambda_m}h_1(s,y)y\,d\nu(y)ds+M\big(h_1\1_{\R_0-\Lambda_m}\big),\\
&\qquad
\sum_{j=1}^kh_2(s_j,y_j)y_j-\int_0^T\int_{\Lambda_m}h_2(s,y)y\,d\nu(y)ds+M\big(h_2\1_{\R_0-\Lambda_m}\big)\Big)\,
 ds_1\cdots ds_kd\nu(y_1)\cdots d\nu(y_k)\Big],
\end{align*}
where  $S_k(\Lambda_m)$
the {\it simplex}
 $$\big\{(t_1,x_1;\dots;t_k,x_k)\in ([0,T]\times \Lambda_m)^k:\ t_1<\dots <t_k\big\}.$$
Therefore, using again that $f^2$ is non negative, $P_{\R_0-\Lambda_m}$--a.s., the
function from $S_k(\Lambda_m)$ into $\R$,
\begin{align*}
(s_1,y_1;\dots;s_k,y_k)\mapsto
f&\big(\sum_{j=1}^kh_1(s_j,y_j)y_j-\int_0^T\int_{\Lambda_m}h_1(s,y)y\,d\nu(y)ds+M\big(h_1\1_{\R_0-\Lambda_m}\big),\\
&\sum_{j=1}^kh_2(s_j,y_j)y_j-\int_0^T\int_{\Lambda_m}h_2(s,y)y\,d\nu(y)ds+M\big(h_2\1_{\R_0-\Lambda_m}\big)\Big)
\end{align*}
is zero a.e. $ds_1\dots ds_k\,\nu(dy_1)\dots \nu(d y_k)$. Hence, by the continuity of
both $f$ and  $h_i(\cdot,y)$, it is deduced that,
fixed   $(y_1,\dots,y_k)$, $\nu^k$--a.e., the
derivatives with respect to the time variables $s$ are also zero.
and then, for every $\ell=1,\dots,k$
\begin{align*}
\sum_{i=1}^2 \partial_i f&\big(\sum_{j=1}^kh_1(s_j,y_j)y_j-\int_0^T\int_{\Lambda_m}h_1(s,y)y\,d\nu(y)+M\big(h_1\1_{\R_0-\Lambda_m}\big),\\
&\sum_{j=1}^kh_2(s_j,y_j)y_j-\int_0^T\int_{\Lambda_m}h_2(s,y)y\,d\nu(y)+M\big(h_2\1_{\R_0-\Lambda_m}\big)\Big)\partial_{s}h_i(s_\ell,y_\ell)y_\ell=0,
\end{align*}
a.e. $ds_1\dots ds_k\,\nu(dy_1)\dots \nu(d y_k)$. Multiplying this quantity by
$k(s_\ell,y_\ell)\big(\frac{s_\ell}{T}-\1_{(t,T]}(s_\ell)\big)$ and
adding with respect to $\ell$, we arrive to
$$\sum_{i=1}^2 \partial_i f\big(M(h_1),M(h_2)\big)\int_{[0,T]\times \Lambda_m}k(s,y)\partial_sh_i(s,y)y\big(\frac{s}{T}-\1_{(t,T]}(s)\big)\, dN(s,y)=0,\ \text{a.s.}$$
Now, thanks to Remark \ref{re:der}, we can  apply the dominated convergence Theorem
to pass to  the limit that expression when $m\to\infty$, and we
obtain that (**) is 0. $\qquad \square$

\bigskip

\bigskip

\begin{proposition}
\label{densitat} ${\cal S}$ is dense in $L^2(\Omega)$.
\end{proposition}

\smallskip

\noindent{\it Proof.} Since the finite linear combinations of
multiple integrals are dense in $L^2(\Omega)$, it is enough to prove
that the random variables of the form $M(E_1)\cdots M(E_n)$ for
$E_1,\dots, E_n\in{\cal B}\big( [0,T]\times \R\big)$  disjoint sets,
with finite $\mu$ measure, can be approximated by elements of $\cal
S$. We proceed in four steps.

\medskip

\noindent{\bf 1.} It suffices to prove the approximation for
  disjoint closed sets. This is proved in the
following way. The measure $\mu$ is $\sigma$--finite in ${\cal
B}\big( [0,T]\times \R\big)$, and hence it is regular in the sense
that for every $A\in{\cal B}\big( [0,T]\times \R\big)$,
$$\mu(A)=\sup\{\mu(C):\, C\subset A,\, C\  \text{closed}\}=
\inf\{\mu(O):\, A\subset O,\, O\  \text{open}\}.$$ So, given
$E_1,\dots,E_n\in {\cal B}\big( [0,T]\times \R\big)$, for every
$\eps>0$, there are closed sets $C_1,\dots, C_n$, with $C_j\subset
E_j$ and
$$\mu(E_j-C_j)<\eps, \ j=1,\dots,n.$$
%From the formula
%$$\prod_{i=1}^na_i-\prod_{i=1}^nb_i=\sum_{i=1}^na_1\cdots a_{i-1}(a_i-b_i)b_{i+1}\cdots b_n,$$
From  the fact that $E_1,\dots,E_n$ are disjoints and $C_j\subset
E_j$, it follows that
$$E\Big[\big(M(E_1)\cdots M(E_n)-M(C_1)\cdots M(C_n)\big)^2\Big]\le n\max_{j=1,\dots, n}\{\mu(E_j)\}^{n-1}
\max_{j=1,\dots, n}\{\mu(E_j-C_j)\},$$ and the claim of this step
follows.
%where $K$ depends only on $n$ and on $\mu(E_1),\dots,\mu(E_n)$.

\medskip

\noindent{\bf 2.} We can assume that $E_1,\dots, E_n$ are disjoint
open sets: From previous step, we can take $E_1,\dots, E_n$ as
disjoint closed sets. By the separability properties of $\R^n$, they
can be separated by open sets, that is, there are disjoint open sets
$O_1,\dots, O_n$ such that $E_j\subset O_j$. By the regularity of
$\mu$, and the same reasoning as before, the property is proved.

\medskip

\noindent{\bf 3.} Let $E_1,\dots, E_n$ be disjoint open sets. Since
every open set on $[0,T]\times \R$ is a countable union of open
rectangles, we can assume that each $E_j$ is a finite union of open
rectangles (with finite $\mu$-measure) of the form $(t',t'')\times
(x',x'')$, or $[0,t)\times (x',x'')$ or $(t,T]\times (x',x'')$ ;
moreover, we can  reduce slightly the size of all the $t$-intervals
in such a way that the measures $\mu(E_j), \ j=1,\dots, n$, do not
change significantly: specifically,    given $\eps>0$, there is
$\delta>0$ such that if we denote $E_j^\delta$ the union with
rectangles $(t'+\delta,t''-\delta)\times (x',x'')$ instead of the
corresponding $(t',t'')\times (x',x'')$ of $E_j$, and similarly for
the other types of rectangles, then
  $$\mu(E_j-E^\delta_j)\le \eps,\ j=1,\dots, n.$$
Let $m_r\in {\cal C}^\infty_0(\R),\ \ r>0$ a family of mollifiers
(see, for example, Adams \cite[Page 29]{adam-75}). For a rectangle
$(t'+\delta,t''-\delta)\times (x',x'')$ write
$$f_r(t,x)=\Big(m_r*\1_{(t'+\delta,t''-\delta)}\Big)(t)\,\1 _{(x',x'')}(x),$$
where the symbol $*$ means convolution.
%\begin{enumerate}[(i)]
%\item For $r<\delta$, $\text{supp}f_r\subset  (t',t'')\times (x',x'')$.
%\item For every $x$, and $r<\delta$, $f_r(\cdot,x)$ is ${\cal C}^\infty((0,T))$ with compact support.
%\item $$\int_{(0,T)\times \R}\big( \partial_t f_r(t,x)\big)^2\,d\mu=
%(t''-t')\int_{(x',x'')} x^2\, d\nu(x) <\infty.$$
%\item By dominated convergence Theorem,
% $$\lim_{r\downarrow 0}\int_{[0,T]\times \R} \Big(f_r(t,x)-
%\1_{(t'+\delta,t''-\delta)\times (x',x'')}(t,x)\Big)^2\, d\mu=0.$$
%\end{enumerate}
Collecting all the $f_r$ corresponding to the same $E^{\delta}_j$,
and using the fact that the support are included in $E_j$ for
$r<\delta$, we obtain that a random variable of the form
$M(g_1)\cdots M(g_n)$, with $g_1,\dots g_n$ with the required
properties, and moreover, with disjoint support, approximates the
product $M(E_1)\cdots M(E_n)$ in $L^2(\Omega).$

\medskip

\noindent{\bf 4.} Finally, use an approximation by ${\cal
C}^\infty_0(\mathbb{R}^n)$ functions to the product function to get the
Proposition. $\qquad \square$

\begin{remark}
\label{suport}
 In the preceding proof we have proved that instead of ${\cal S}$ we can use the sum of the smooth random variables of the form
$f(M(h_1),\dots,M(h_n))$ with $h_1,\dots, h_n$ with disjoint support
to define the local Malliavin derivative.
\end{remark}

\bigskip

Next proposition gives a Malliavin  integration by parts formula:

\begin{proposition}\label{pro:duality}
Let $F\in \mathcal{S}$, and $g$ be a measurable and bounded function
on $[0,T].$ Then
\begin{align}
\label{parts} \E\Big[&\int_0^T ( D^{\Lambda}_{t} F )g(t)
dt\Big]\\
&=\E\Big[F\Big(\1_{\Lambda}(0)\sigma \int_0^T
g(s)dW_s+\int_{[0,T]\times(\Lambda\cap\R_0)}
\big(g(s)-\frac{1}{T}\int_0^T
g(t)dt\big)k(s,y) dN(s,y)\Big)\Big]\\
&\qquad\qquad-\E\Big[F\Big(
\int_{[0,T]\times(\Lambda\cap\R_0)}\partial_s k(s,y) \big[\int_0^T
g(t)(\frac{s}{T}-{\1}_{[0,s]}(t))dt\big]dN(s,y)\Big)\Big].
\end{align}
\end{proposition}
\textit{Proof:}

First we consider the case  $F=f(M(h)),$ where $f$ and $h$ are as in
(\ref{smoothfunctional}) and
 assume  $\nu(\Lambda)<\infty.$ Then,
\begin{eqnarray}
\label{domini:1}\lefteqn{
\E\left[\int_0^T  (D^{\Lambda}_{t} f(M(h))) g(t) dt\right]}\nonumber\\
&=&\1_{\Lambda}(0) \E \left[ \sigma\int_0^T f'(M(h))  h(t,0)
g(t) dt\right]\nonumber\\
&& +\E\left[\int_0^T f'(M(h))
\left(\int_{[0,T]\times(\Lambda\cap\R_0)} k(s,y)\partial_s h(s,y) y
(\frac{s}{T}-\1_{[t,T]}(s))dN(s,y)\right)g(t)
 dt\right].
 \end{eqnarray}

From the integration by parts formula in the Wiener case,
$$\E\left[\int_0^T D_t^{\{0\}}f(M(h))g(t)dt\right]=
\E\left[ f(M(h))M(g(t)\otimes\1_{\{0\}}(x))\right]$$ and the
 first term  of the right hand side of (\ref{domini:1}) is equal to
$$\1_{\Lambda}(0) \E[f(M(h)) M( g \otimes\1_{\Lambda \cap \{0\}})].$$

Now, using Fubini theorem, the second term  can be written as
\begin{align*}
\E\Big[f'(M(h))&\Big( \int_{[0,T]\times(\Lambda\cap {\mathbb
R}_0)}\Big\{ k(s,y)\partial_s h(s,y)y\big(\frac{s}{T}\int_0^T g(t)dt
- \int_0^s g(t) dt\big)\Big\} dN(s,y)\Big)\Big],
\end{align*}
which, by the independence of $W$, $N^{\Lambda}$ and $N^{{\mathbb
R}_0-\Lambda}$ (using the same notations as in the previous proof),
and the symmetry of the function $\sum_{j=1}^nh(s_j,y_j)y_j$,  is
equal to
\begin{align*}
\sum_{n=1}^{\infty}& \Prob\{N^{\Lambda\cap\R_0}_T=n\}
\frac{1}{\nu(\Lambda)^nT^n}\\
&\times \int_{[0,T]^n} \int_{(\Lambda\cap\R_0)^n}\Big\{
 \E_{W,N^{\R_0-\Lambda}}\Big[f'(\sum_{j=1}^nh(s_j,y_j)y_j
-\int_0^T\int_\Lambda
h(r,z)z\,dr d\nu(z)\\
&\hspace{1cm}+\1_{\Lambda}(0) M(h\1_{\{0\}})
+ M(h\1_{\R_0-\Lambda}))\Big]\\
&\times  \sum_{j=1}^n k(s_j,y_j)\partial_s
h(s_j,y_j)y_j\Big(\frac{s_j}{T}\int_0^T g(t)dt- \int_0^{s_j} g(t)
dt\Big) \Big\}\, ds_1\cdots ds_n d\nu(y_1)\cdots d\nu(y_n).
\end{align*}
Integrating by parts each integral
\begin{align*}
\int_0^T \E_{W,N^{\R_0-\Lambda}}&\Big[f'(\sum_{j=1}^nh(s_j,y_j)y_j
-\int_0^T\int_\Lambda h(r,z)z\,dr d\nu(z)
+\cdots)\Big]\\
&\times k(s_i,y_i)\partial_s h(s_i,y_i)y_i\Big(\frac{s_i}{T}\int_0^T
g(t)dt- \int_0^{s_i} g(t) dt\Big)ds_i,
\end{align*}
and reordering terms, we arrive to
\begin{align*}
\E\Big[f(M(h)) & \int_{[0,T]\times
(\Lambda\cap\R_0)}k(s,y)\Big(g(s)-\frac{1}{T}\int_0^T g(t)dt\Big)
dN(s,y)\Big]
\\
&-\E\Big[f(M(h)) \int_{[0,T]\times (\Lambda\cap\R_0)}\partial_s
k(s,y)\Big(\int_0^T g(t)(\frac{s}{T}-{\1}_{[0,s]}(t))dt\Big)
dN(s,y)\Big],
\end{align*}
which implies that the result holds in the case that
$\nu(\Lambda)<\infty$.

To deal with the case $\nu(\Lambda)=\infty$,  let
${\Lambda}_n=\{0\}\cup\{\frac1n<|x|<n\}$. Then, changing $\Lambda$ by
${\Lambda}_n\cap \Lambda,$ and letting  $n\rightarrow \infty,$ we obtain
the result.

Now, consider  $F=f(M(h_1),\dots,M(h_n))$, with $h_1,\dots,h_n$ with
disjoint support. Since the random variables $M(h_1),\dots, M(h_n)$
are independent, we can compute the expectations in (\ref{parts})
iteratively, and the result follows. Finally, by linearity, the
equality  (\ref{parts}) is extended to a sum of smooth random
variables of the form $f(M(h_1),\dots,M(h_n))$, with $h_1,\dots,h_n$
with disjoint support, which generate ${\cal S}$, in agreement with
Remark \ref{suport}. $\qquad\square$

\bigskip

\bigskip

The following result is a consequence of Proposition
\ref{pro:duality} applied to $FG$ and of the fact that $D^{\Lambda}$
is a true derivative operator, which follows from Definition
\ref{localderivative}.
\begin{corollary} \label{cor:prod}
Let $g$ measurable and bounded, and $F,G\in {\cal S}$. Then

\begin{eqnarray*}
\lefteqn{\E\left(G\int_0^T ( D^{\Lambda}_{t} F )g(t) dt+F\int_0^T (
D^{\Lambda}_{t} G )g(t)
dt\right)}\\
&&=\E\left[FG\left(\1_{\Lambda}(0) \int_0^T \sigma g(s)dW_s+\
\int_{[0,T]\times (\Lambda\cap\R_0)} \left(g(s)-\frac{1}{T}\int_0^T
g(t)dt\right)k(s,y) d
N(s,y)\right.\right.\\
&&\left. \left.-\int_{[0,T]\times (\Lambda\cap\R_0)}\partial_s
k(s,y) \left(\int_0^T
g(t)(\frac{s}{T}-{\1}_{[0,s]}(t))dt\right)dN(s,y)\right) \right].
\end{eqnarray*}
\end{corollary}

\bigskip

\bigskip

Now we can proceed as in Nualart \cite{Nu-95} (pag. 26) to show  the
following result, taking into account that the bounded functions are
dense in $L^2 ([0,T]).$

\begin{corollary} The operator
 $D^{\Lambda},$ from $L^2
(\Omega)$ into $L^2 ([0,T]\times \Omega),$
is unbounded, densely defined and closable.
\end{corollary}

In particular, we have that the operator $D^{\Lambda}$ has a closed
extension, which is also written   by  $D^{\Lambda}$. The domain of
this operator is denoted by $\text{Dom}\,D^{\Lambda}$. Note that $F$
belongs to  $\text{Dom}\,D^{\Lambda}$ if and only if there is a
sequence $\{F_n: n\ge 1\}\subset{\cal S}$ such that $F_n$ converges
to $F$ in $L^2(\Omega)$, when $n\to\infty$, and $\{D^{\Lambda}F_n:
n\ge 1\}$ converges to some process $Y$ in $L^2
([0,T]\times\Omega).$ In this case, we have that $D^{\Lambda}F=Y$.

\bigskip

The following  rules of derivation are well known in Brownian
Malliavin calculus; they can be translated to our context with the
same proof.

\begin{proposition}
\null

\noindent{\bf 1.} ({\it Chain rule}). Let $f:\R^n\to\R$ be a
continuously differentiable function with bounded partial
derivatives, and let $\bs F=(F_1,\dots,F_n)$ a random vector such
that $F_j\in {\rm Dom}\, D^\Lambda$, for $j=1,\dots, n$. Then $f(\bs
F)\in {\rm Dom}\, D^\Lambda$ and
$$D^\Lambda\big(f(\bs F)\big)=\sum_{j=1}^n\partial_j f(\bs F) D^\Lambda F_j.$$

\noindent{\bf 2.} Let $F,G\in {\rm Dom}\, D^\Lambda$ such that $ G\,
D^\Lambda F$ and $F\,  D^\Lambda G$ belong to $L^2([0,T]\times \Omega)$. Then
$F\,G\in {\rm Dom}\, D^\Lambda$ and
$$ D^\Lambda(F\,G)=G\,  D^\Lambda F+ F\,  D^\Lambda G.$$

\end{proposition}

\subsection{Absolute continuity via the operators $\bs{D^{\Lambda}}$}
The chain rule allows to prove a criterion  for the absolutely
continuity of a functional $F$ in a similar way that in the Brownian
case. Specifically,  the result obtained by  Carlen and Pardoux
\cite[Theorem 4.1]{Ca-Pa-90}  for the Poisson process can be
extended to the Lévy case; the proof  is the same as Carlen and
Pardoux,   and omitted.
\begin{proposition}\label{pro:abcon}
Let $\Lambda\in{\cal B}(\R)$, $k\in{\cal K}$ and $F\in{\rm Dom}\,
D^\Lambda$ such that
$$\int_0^T(D_t^{\Lambda} F)^2 dt>0$$
a.s. on a measurable set $A\in{\cal F}$. Then, $\Prob\circ F^{-1}$
is absolutely continuous on $A,$ that is,  $\lambda(B)=0$ implies
$\Prob(\{F\in B\}\cap A)=0$.
\end{proposition}
It is worth to remark that the criterion is true for every set
$\Lambda$ and weight function $k\in {\cal K}$. This is very
interesting  for applications because we can choose an appropriate
$\Lambda$ and $k$ depending of $F$.   It will be used in Section
\ref{sec:4} to see that the solutions of some stochastic differential
equations have a density on a measurable set.

\section{Multiple integrals over sets with finite L\'evy measure and its derivatives}
\label{sec:multiple} In this section we study   multiple integrals
defined on a set with finite Lévy measure. Such integrals can be
computed pathwise, and its properties can be proved by combinatorial
methods.

 Let $\Theta\in{\cal B}(\R)$ be a bounded set such that $0\not\in \overline \Theta$ (closure of $\Theta$); then $\nu(\Theta)<\infty$.
As we have done  before, consider
$$N^\Theta(B)=\#\{t:\ (t,\Delta X_t)\in B \ \text{and} \ \Delta X_t\in \Theta\},  \ \text{for}\  B\in {\cal B}\big((0,\infty)\times \R_0 \big),$$
that is a Poisson random measure with intensity $\lambda\otimes
\nu_\Theta$, where $\nu_\Theta(C)=\nu(\Theta\cap C)$, for $C\in{\cal
B}(\R).$ Write
 $$N^\Theta_T=N^\Theta\big([0,T]\times\Theta\big) <\infty, \ \text{a.s.}$$
% Since in the interval $[0,T]$ a L\'evy process can has only a finite number of jumps of size $\vert \Delta X_t\vert>\epsilon$ for $\eps>0$ arbitrary,
%$$N^\Theta_T<\infty, \ \text{a.s.}$$
Then we can order the jumps in the interval $[0,T]$, say
$T_1<\cdots<T_{N^\Theta_T}$, and write
$$N^\Theta=\sum_{j=1}^{N^\Theta_T}\delta_{(T_j, \Delta X_{T_j})}.$$
Let $\phi:[0,T]\times \R_0\to\R$ be a measurable function (indeed,
it is only needed that $\phi$ is defined in $[0,T]\times\Theta$) integrable with respect to $N^\Theta$ and consider
 $$J_1^\Theta(\phi):=\int_{[0,T]\times \mathbb{R}_0} \phi(t,x)\, dN^\Theta(t,x).$$
This integral has a compound Poisson distribution with
characteristic function (see, for example, Sato \cite[Proposition
19.5]{Sa-99})
$$\E\big[e^{iuJ_1^\Theta(\phi)}\big]=\exp\big\{\int_{[0,T]\times\Theta}\big(e^{iu\phi(t,x)}-1\big)\,dtd\nu (x)\big\}.$$
Thus, if $\phi\in L^p([0,T]\times \R_0,\lambda\otimes\nu_\Theta)$,
then $J_1^\Theta(\phi)\in L^p(\Omega)$, for $p\ge 1$.

% If $h\in L^2([0,T]\times R_0,\lambda\otimes\nu_\Theta)$,
%then $I:=\int_{[0,T]\times R_0} h(s,x)\, dN^\Theta(s,x)\in L^2(\Omega)$ and
%$$\E[I]=\int_0^T\int_\Theta{\Theta} h(s,x)\,dsd\nu(x)
%\quad \text{and}\quad \E[I^2]=\int_0^T\int_{\Theta} h^2(s,x)\,dsd\nu(x)+\big(\int_0^T\int_{\Theta} h(s,x)\,dsd\nu(x)\big)^2.$$

\bigskip

Also multiple integrals can be considered. Recall the notations introduced
in the proof of Proposition \ref{definit}: let  $S_n(\Theta)$
be the {\it simplex}
 $$\big\{(t_1,x_1;\dots;t_n,x_n)\in ([0,T]\times \Theta)^n:\ 0\le t_1<\dots <t_n\le T\big\}.$$
 For $\phi: S_n(\Theta)\to\R$, define
\begin{align*}
 J_n^\Theta(\phi)&=\int\cdots\int_{S_n(\Theta)} \phi (t_1,x_1;\dots;t_n,x_n)\, dN^\Theta(t_1,x_1)\cdots d N^\Theta(t_n,x_n)\\
& =\sum_{1\le i_1<\cdots<i_n\le N^\Theta_T}\phi(T_{i_1},\Delta
X_{T_{i_1}};\dots;T_{i_n},\Delta X_{T_{i_n}})\, \1_{\{N^\Theta_T\ge
n\}}.
\end{align*}
%with the convention that the sum is zero if $N^\Theta_T<n.$

\begin{proposition} For every $p\ge 1$, if   $\phi\in L^p\big(S_n(\Theta),(\lambda\otimes \nu_\Theta)^n\big)$,
then $ J_n^\Theta(\phi)\in L^p(\Omega)$, and
\begin{equation}
\label{fita} \E\big[\big\vert J_n^\Theta(\phi)\big\vert^p\big]\le
C_{p,n} \int_{S_n(\Theta)}\vert \phi(t_1,x_1;\dots;t_n,x_n)\vert^p\,
dt_1 \cdots dt_n d\nu(x_1)\cdots d\nu(x_n),
\end{equation}
where the constant $C_{n,p}$ depends only on $n$ and $p$.

\end{proposition}

\medskip

\noindent{\it Proof.} The main tool of the proof is that the  jump times
 of a compound Poisson
 process follow a Dirichlet distribution, independent of the size of jumps,
 which are also independent
 and each of one has  law $P\{\Delta X_{T_j}\in A\}=\nu(A)/\nu(\Theta)$, for $A\in{\cal B}(\Theta).$
To simplify the notations, we assume $T=1$. By Jensen inequality,
$$\vert J_n^\Theta(\phi)\big\vert^p\le (\Big(N^{(\Theta)}\otimes
\cdots\otimes N^{\Theta})(S_n(\Theta))\Big)^{p-1} \int_{S_n(\Theta)} \vert \phi
(t_1,x_1;\dots;t_n,x_n)\vert^p\, dN^\Theta(t_1,x_1)\cdots d
N^\Theta(t_n,x_n).$$ Hence,
\begin{align*}
\E\big[&\big\vert J_n^\Theta(\phi)\big\vert^p\big]=\sum_{k\ge n}
\E\Big[\big\vert J_n^\Theta(\phi)\big\vert^p\,\vert\,N_1^\Theta=k\Big]\,\Prob\{N_1^\Theta=k\}\\
&\le e^{-\nu(\Theta)}\sum_{k\ge n}\frac{\nu(\Theta)^k}{k!}\,k^{p-1}
\frac{1}{\nu(\Theta)^n}\,
\E\Big[\sum_{1\le {i_1}<\cdots<{i_n}\le k}\big\vert \phi(T_{i_1},\Delta X_{T_{i_1}};\dots;T_{i_n},\Delta X_{T_{i_n}})\Big\vert^p\Big]\\
&= e^{-\nu(\Theta)}\sum_{k\ge n}\frac{\nu(\Theta)^k}{k!}\,k^{p-1}
\frac{1}{\nu(\Theta)^n}\sum_{1\le {i_1}<\cdots<{i_n}\le k}
\int_{\Theta^n}\E\Big[\big\vert
\phi(U_{(i_1)},x_1;\dots;U_{(i_n)},x_n)\big\vert^p\Big]
 d\nu(x_1)\dots d\nu(x_n),
\end{align*}
where   $(U_{(i_1)},\dots,U_{(i_n)})$ is the ordered
$(i_1,\dots,i_n)$--statistic from a  sample
 $U_{1},\dots,U_{k}$ of a uniform law in $(0,1)$.
%; such density is
%\begin{align*}
%f_{i_1,\dots,i_n;k}(t_1,\dots,t_n)&=\frac{k!}{(i_1-1)!(i_2-i_1-1)!\cdots (i_n-i_{n-1}-1)(k-i_{n})!}\\
%&\qquad \cdot t_1^{i_1-1}(t_2-t_1)^{i_2-i_1-1}\cdots (t_n-t_{n-1})^{i_n-i_{n-1}-1}(1-t_n)^{k-i_{n}},
%\end{align*}
%for $0<t_1<\cdots<t_n<1.$
Extend $\phi$ to the whole $(0,1)^n$ by symmetry, and denote by
${\cal C} _n^k$  the family of all subsets of $\{1,\dots,k\}$ of
size $n$. Fixed $x_1,\dots, x_n\in \Theta,$
\begin{align*}
\E\Big[\sum_{1\le {i_1}<\cdots<{i_n}\le k} & \big\vert
\phi(U_{(i_1)},x_1;\dots ;U_{(i_n)},x_n)\big\vert^p\Big]
=\E\Big[\sum_{\{i_1,\dots,i_n\}\in {\cal C} _n^k} \big\vert
\phi(U_{i_1},x_1;
\dots;U_{i_n},x_n)\big\vert^p\Big]\\
& =\binom{k}{n}\int_{(0,1)^n}\big\vert \phi(t_1,x_1,\dots,t_n,x_n)\big\vert^p\, dt_1\cdots dt_n\\
&= \frac{k!}{(k-n)!}\int_{0<t_1\cdots<t_n<1}\big\vert
\phi(t_1,x_1,\dots,t_n,x_n)\big\vert^p\, dt_1\cdots dt_n.
\end{align*}
%where ${\cal C} _n^k$ is the family of all subsets of $\{1,\dots,k\}$ of size $n$.
%Hence,
%\begin{align*}
%\E\Big[&\sum_{1\le {i_1}<\cdots<{i_n}\le k}
%\big\vert \phi(U_{(i_1)},x_1,\dots,U_{(i_n)},x_n)\big\vert^p\Big]=

%\end{align*}
Finally,
\begin{align*}
\E\big[&\big\vert J_n^\Theta(\phi)\big\vert^p\big] \le
\int_{0<t_1\cdots<t_n<1}\int_{\Theta^n}
\big\vert \phi(t_1,x_1,\dots,t_n,x_n)\big\vert^p\, dt_1\cdots dt_n d\nu(x_1)\dots d\nu(x_n)\\
&\times e^{-\nu(\Theta)}\sum_{k\ge
n}\frac{\nu(\Theta)^{k-n}}{(k-n)!}\,k^{p-1},
\end{align*}
and the result follows. $\quad \square$

\begin{corollary}
\label{uniform} Let $\phi_k,\phi\in
L^p\big(S_n(\Theta),(\lambda\otimes \nu_\Theta)^n\big)$, such that
$$\lim_k\phi_k(t_1,x_1;\dots;t_n,x_n)=\phi(t_1,x_1;\dots;t_n,x_n), \ \text{uniformly}.$$
Then
$$\lim_k J_n(\phi_k)= J_n(\phi), \ \text{in}\ L^p(\Omega).$$
\end{corollary}

\bigskip

 A crucial result in multiple integration theory is the product formula (see Lee and Shih \cite{chinos} and Tudor \cite{tudo-97}). For this particular
  case, the formula is simpler and the proof straightforward:

\begin{proposition}
\label{product}

Consider $\phi_n: S_n(\Theta)\to\R$ and $\phi_1:[0,T]\times
\R_0\to\R$. Then,
\begin{equation}
\label{form-prod}
J_n^\Theta(\phi_n)\,J_1^\Theta(\phi_1)=J_{n+1}^\Theta(\phi_n\wt\otimes
\phi_1)+J_n^\Theta(\phi_n*\phi_1),
\end{equation}
where
$$\phi_n\wt\otimes \phi_1 (t_1,x_1;\dots;t_{n+1},x_{n+1})=\sum_{j=1}^{n+1}\phi_n(t_1,x_1;\dots;\wh{t_j,x_j};\dots ;t_{n+1},x_{n+1})\phi_1(t_j,x_j),$$
where $\wh{t_j,x_j}$ means that this pair is missing, and
$$\phi_n* \phi_1 (t_1,x_1;\dots;t_{n},x_{n})=\sum_{j=1}^{n}\phi_n(t_1,x_1;\dots ;t_{n},x_{n})\phi_1(t_j,x_j).$$
\end{proposition}

\medskip

\bigskip
\begin{remark}
\label{symm} Assume that  $\phi$ is defined in all $([0,T]\times
\Theta)^n$.  The function $\phi_n\wt\otimes \phi_1$ is not the
standard symmetrization of $\phi_n\otimes\phi_1$, which is given by
$$\wt{\phi_n\otimes \phi_1} (t_1,x_1;\dots;t_{n+1},x_{n+1})=
\frac{1}{(n+1)!}\sum_{\sigma}\phi_n(t_{\sigma_1},x_{\sigma_1};\dots;t_{\sigma_n},x_{\sigma_n})\phi_1(t_{\sigma_{n+1}},x_{\sigma_{n+1}}),$$
where the sum ranges over the set of permutations of $n+1$ elements.
We have the relation
$$\wt{\phi_n\otimes \phi_1}=\frac{1}{n+1}\,\wt\phi_n\wt\otimes \phi_1,\ \text{where}\ \wt\phi_n \ \text{is the symmetrization of}\ \phi_n.$$

\end{remark}

\begin{lemma}
\label{lema-mult} Let $\Theta\in{\cal B}(\R)$ be a bounded set such
that $0\not\in \overline \Theta$, and $\Lambda \in{\cal B}(\R)$ with
$\Theta\subset \Lambda$. Consider $f\in {\cal C}^\infty_p(\R).$ Then
$ f(N^\Theta_T)\in{\rm Dom}\, D^{\Lambda}$ and
$$D^\Lambda f(N^\Theta_T)=0.$$
\end{lemma}

\medskip

\noindent{\it Proof.} Write
\begin{align*}
N^\Theta_T&=\int_{[0,T]\times \R_0} dN^\Theta(t,x)=\int_{[0,T]\times \R_0} \frac{1}{x}\1_{\Theta}(x)x\,dN(t,x)\\
&=\int_{[0,T]\times \R_0} \frac{1}{x}\1_{\Theta}(x)x \,d\wt
N(t,x)+\int_0^T\int_\Theta\,dtd\nu(x).
\end{align*}
Hence, $N^\Theta_T\in {\cal S}$, and $D^\Lambda N^\Theta_T=0$. Thus
also $ f(N^\Theta_T)\in{\cal S}$ and $D^\Lambda
f(N^\Theta_T)=0.\qquad \square$
\bigskip

Given a set  $A\subset \N$, the random variable $\1_{\{N^\Theta_T\in
A\}}$ can be conveniently approximated by functionals
$f(N_t^\Theta)\in {\rm Dom}\, D^{\Lambda}$. From the closability of
$D^\Lambda$ we get the following corollary:
\begin{corollary}
\label{indicador} For every set $A\subset \N$, $\1_{\{N^\Theta_T\in
A\}}\in {\rm Dom}\, D^{\Lambda}$ and $D^\Lambda \1_{\{N^\Theta_T\in
A\}}=0.$
\end{corollary}

\begin{theorem}
\label{theo-mult} Let $\Theta$ and $\Lambda$ as in  Lemma
\ref{lema-mult}, and $\phi\in L^2\big(S_n(\Theta),
(\lambda\otimes\nu)^n\big)$  such that for every $(x_1,\dots,x_n)\in
\Theta^n$, $\phi(\cdot, x_1;\dots;\cdot,x_n)\in {\cal
C}^\infty\big(\overline{S_n}\big),$ where $\overline{S_n}=\{0\le
t_1\le \cdots\le t_n\le T\}$.
 Then  $J_n^\Theta(\phi)\in {\rm Dom}\, D^{\Lambda}$ and
\begin{equation}
\label{derivada-int}
D^{\Lambda}_tJ_n^\Theta(\phi)=\sum_{j=1}^nJ_n^\Theta\Big(k(s_j,x_j)\partial_{s_j}\phi(s_1,
x_1;\dots;s_n,x_n)\big(\frac{s_j}{T}-\1_{(t,T]}(s_j)\big)\Big).
\end{equation}
\end{theorem}

\medskip

\noindent{\it Proof.} We proceed by induction. For $n=1$, as in the
previous lemma,  write
\begin{align*}
J_1^\Lambda(\phi)&=\int_{[0,T]\times \R_0}\phi(s,x)\, d N^\Theta(s,x)\\% \,\1_{\{N^\Theta_T\ge 1\}}
&=M(\1_\Theta(x)\phi(s,x)/x)+\int_0^T\int_\Theta \phi(s,x)\, dtd \nu(x). %\, \1_{\{N^\Theta_T\ge 1\}}.
\end{align*}
Its derivative is
$$\int_{[0,T]\times \mathbb{R}_0}\1_\Theta(x)k(s,x)\partial_s\phi(s,x)\big(\frac{s}{T}-\1_{(t,T]}(s)\big)\, dN(s,x),$$
%J_1^\Theta(k(s,x)\partial_s\phi(s,x)\big(\frac{s}{T}-\1_{(t,T]}(s)).$$
which coincides with (\ref{derivada-int}).

Now, consider $\phi_n \in{\cal C}_0^\infty\big(S_n(\Theta)\big)$ and
$\phi_1\in {\cal C}_0^\infty\big((0,T)\times R_0\big)$. Then the
result is obtained using
 the product formula (\ref{form-prod})
\begin{equation}
\label{prod-apl} J_{n+1}^\Theta(\phi_n\wt\otimes
\phi_1)=J_n^\Theta(\phi_n)\,J_1^\Theta(\phi_1)-J_n^\Theta(\phi_n*\phi_1),
\end{equation}
checking that all terms in the right hand side are in ${\rm Dom}
D^\Lambda$ and using the induction hypothesis.

For a general $\phi \in L^2\big(S_{n+1}(\Theta),
(\lambda\otimes\nu)^{n+1}\big)$ with the conditions  of the theorem,
first extend this
 function by symmetry to  $([0,T]\times \Theta)^n$
Since  $\Theta$  is a bounded set by hypothesis, using a mollifier
family $m_\eps\in {\cal C}_0^\infty(\R^{n+1})$,  define
\begin{align*}
\phi_\eps(t_1,x_1;\dots;t_{n+1},x_{n+1})&=\int_{\R^{n+1}}\phi(t_1,y_1;\dots;t_{n+1},y_{n+1})m_\eps(x_1-y_1,\dots,x_{n+1}-y_{n+1})\\
&\qquad \qquad d\nu(y_1)\cdots d\nu(y_{n+1}).
\end{align*}
Then,   $\phi_\eps \in {\cal C}_0^\infty(\R^{2n+2}),$
$\lim_{\eps\downarrow 0}\phi_\eps=\phi,\ \text{uniformly}\quad\text
{and}\quad \lim_{\eps\downarrow
0}\partial_{t_j}\phi_\eps=\partial_{t_j}\phi,\ \text{uniformly}.$
Hence, by Corollary \ref{uniform}, $\lim_{\eps\downarrow
0}J_{n+1}(\phi_\eps)=J_{n+1}(\phi)$ in $L^2(\Omega)$, and also the
right hand side of (\ref{derivada-int}) corresponding to $\phi_\eps$
converges in $L^2([0,T]\times \Omega)$ to the corresponding to
$\phi$. Thus, it suffices to prove the formula (\ref{derivada-int})
for $\phi\in {\cal C}_0^\infty(\R^{2n+2})$. For such $\phi$ we can
apply Theorem 1, page 65 of Yosida \cite{Yosida}, that says that we
can approximate $\phi$  by finite sums
$\sum_{j=1}^k\phi_n^{(j)}\phi_1^{(j)}$, where $\phi_n^{(j)}\in {\cal
C}_0^\infty(\R^{2n})$ and  $\phi_1^{(j)}\in {\cal
C}_0^\infty(\R^{2})$, and also the derivatives of the sums converge
uniformly to the derivatives of $\phi$. Indeed, that theorem by
Yosida \cite{Yosida} shows  the result for complex valued functions,
however, it can be applied to real values functions and, in that
case,  the approximating functions are also real. From Corollary
\ref{uniform}, the symmetrization property of Remark \ref{symm} and
the previous induction step, the proof is complete. $\qquad \square$

\begin{corollary}\label{cor:derjt}
Let $\Theta$, $\Lambda$ and $\phi$  be as in Theorem
\ref{theo-mult}. With the previous notations,
 $$\phi(T_1,\Delta
X_{T_1};\dots,T_n,\Delta X_{T_n})\1_{\{N^\Theta_T\ge n\}}\in{\rm
Dom}\,D^{\Lambda}$$
 and, over $\{N^\Theta_T\ge n\}$,
$$D_t^{\Lambda}\phi(T_1,\Delta X_{T_1};\dots;T_n,\Delta X_{T_n})
=\sum_{j=1}^n k(T_j,\Delta X_{T_j})\partial_j\phi(T_1,\Delta
X_{T_1};\dots,T_n,\Delta X_{T_n})
\big(\frac{T_j}{T}-\1_{[0,T_j]}(t)\Big).$$
\end{corollary}

 \noindent\textit{Proof:} We only need to
observe that
$$\phi(T_1,\Delta
X_{T_1};\dots,T_n,\Delta X_{T_n})\1_{\{N^\Lambda_T\ge
n\}}=\sum_{\ell=n}^\infty
J_\ell(\psi_\ell)\1_{\{N_T^{\Lambda}=\ell\}},$$ where
$$\psi_\ell(s_1,y_1;\ldots;s_\ell,y_\ell)=\phi(s_1,y_1;\ldots;s_n,y_n), \quad \ell\ge n,$$
and apply Corollary \ref{indicador} and Theorem \ref{theo-mult}.
$\qquad\square$

\vskip .2cm

\noindent\textbf{Remark.} Note that this result yields that
$D^{\{1\}}$, with $k=1$, agrees with the operator introduced in
Carlen and Pardoux \cite{Ca-Pa-90} for the Poisson process.
Consequently, this paper
 extends the approach given in \cite{Ca-Pa-90} to L\'evy
processes. Also, in general, with the previous hypothesis,  the jump
time $T_j$ is in the domain of $D^{\Lambda}$ and
$$D_t^{\Lambda}T_j=k(T_j,\Delta
X_{T_j})\left(\frac{T_j}{T}- \1_{[0,T_j]}(t)\right).$$

\section{Absolutely continuity for stochastic differential
equations}\label{sec:4}

In this section we use Proposition \ref{pro:abcon} and Theorem
\ref{theo-mult}
 to find  conditions that guarantee that the
solution at time $T$  of  some stochastic differential equations
driven by a Lévy process has  density with respect to the Lebesgue
measure. The key point is that we can choose the convenient set
$\Lambda$ and weight $k(t,x)$ for each type of equation.

\subsection{An equation driven by a L\'evy process with continuous part}
Here we consider the solution of the stochastic differential
equation with an additive jump noise and a Wiener stochastic
integral of the form
$$Z_t=x_0+\int_0^t b(Z_s)ds+\int_0^t\sigma(Z_s)dW_s+\int_{(0,t]\times\R_0}l(y)ydN(s,y),\quad
t\in[0,T],$$ where $l\in{\cal K}$ and  the coefficients $b$ and
$\sigma$ are differentiable on $\R$ with bounded derivatives. Then,
choosing $\Lambda=\{0\}$ and $k\in{\cal K}$, and  proceeding as in
Nualart \cite{Nu-95} (Section 2.3), we can see that
$$D^{\{0\}}_tZ_T=\sigma(Z_t)\exp\left(\int_0^t\sigma'(Z_s)dW_s+\int_0^t\left\{b'(Z_s)-\frac12
(\sigma'(Z_s))^2\right\}ds\right),\quad t\in[0,T].$$ Thus we can
state the following theorem:
% where we give a set in which the
%integral of the square of the last derivative  is positive.
\begin{theorem}
The random variable $Z_T$  is absolutely continuous with respect to
the Lebesgue measure on the set $\{S<T\}$, where
$$S=\inf\left\{t\in[0,T]: \int_0^t \1_{\{\sigma(Z_s)\neq
0\}}ds>0\right\}\wedge T.$$
\end{theorem}
\vskip .2cm \noindent\textbf{Remark} As in \cite{Nu-95}, we can see
that last theorem also holds in the case that the coefficients are
only Lipschitz functions with linear growth.
\subsection{A pure discontinuous equation with a monotone
drift}\label{sec:4.2}
 In this section we consider the following equation with an additive
 jump noise
\begin{equation}\label{eq:mondri}
Z_t=x+\int_0^tf(Z_s)ds+\int_{[0,t]\times\R_0}h(y)dN(s,y),\quad
t\in[0,T],
\end{equation}
where $x\in\R$, $f:\R\rightarrow\R$ is a function with a bounded
derivative and $h\in L^2(\R_0,\nu)\cap L^1(\R_0,\nu)$. It is
well-known that equation (\ref{eq:mondri}) has a unique
square-integrable solution.

\bigskip

In order to calculate $D^{\R_0} Z_T$,  we need the following result.
\begin{proposition}\label{pro:convsol} Let
$\Theta_m=\{x\in \R: \ 1/m<\vert x\vert <m\},$  and
\begin{equation}
\label{eq-approx}
Z_t^{(m)}=x+\int_0^tf(Z_s^{(m)})ds+\int_{[0,t]\times\Theta_m}h(y)dN(s,y),\quad
t\in[0,T].
\end{equation} Then, as $m\rightarrow \infty$, $Z_t^{(m)}$
converges to $Z_t$ in $L^2(\Omega)$,  for every $t\in[0,T]$. That
convergence is also a.s. for every $t\in[0,T]$ a.e.
\end{proposition}

\noindent\textit{Proof:} Let $M>0$ be a bound for the derivative of
the function $f$. Then
$$\E\left(\left|Z_t^{(m)}-Z_t\right|^2\right)\le 2M^2T\int_0^t\E\left(\left|Z_s^{(m)}-
Z_s\right|^2\right)
ds+2\E\left(\left|\int_0^t\int_{[0,t]\times(\R_0-\Theta_m)}h(y)dN(s,y)\right|^2\right).$$
The result is a consequence of Gronwall's lemma and of the fact that
\begin{eqnarray*}
\lefteqn{\E\left(\left|\int_{[0,t]\times(\R_0-\Theta_m)}h(y)dN(s,y)\right|^2\right)}\\
&=&\E\left(\left|\int_{[0,t]\times(\R_0-\Theta_m)}h(y)d{\tilde
N}(s,y)+\int_0^t\int_{\R_0-\Theta_m}h(y)\nu(dy)ds\right|^2\right)\\
&\le&2\int_0^T\int_{\R_0-\Theta_m}h(y)^2\nu(dy)ds +
2\left(\int_0^T\int_{\R_0-\Theta_m}|h(y)|\nu(dy)ds\right)^2
\end{eqnarray*}
goes to zero as $m\rightarrow\infty$ because $\Theta_m\uparrow
\R_0$. Finally the second claim of this result follows using a
similar argument. Thus, the proof is complete.  $\qquad \square$

\bigskip

\bigskip

Now we follow a similar reasoning as Carlen and Pardoux
\cite{Ca-Pa-90}.  Consider the flow $\{\Phi_t(s,x): 0\le s\le t\le T
\ \textrm{and}\ x\in\R\}$ associated with equation
(\ref{eq:mondri}). It means, $\Phi_t(s,x)$ is the unique solution to
the equation
$$\Phi_t(s,x)=x+\int_s^t f(\Phi_u(s,x))du,\quad t\in[s,T].$$
Note that the partial derivatives   of the flow with respect to time
is
$$\partial_t\Phi_t(s,x)=f(\Phi_t(s,x)).$$
The partial derivative with respect the initial condition satisfies
the equation
$$\partial_x\Phi_t(s,x)=1+\int_s^t
f'(\Phi_u(s,x))\partial_x\Phi_u(s,x)du$$ and so
$$\partial_x\Phi_t(s,x)=\exp\left(\int_s^tf'(\Phi_u(s,x))du\right).$$
Similarly we have
$$\partial_s\Phi_t(s,x)=-f(\Phi_s(s,x))+\int_s^t f'(\Phi_u(s,x))\partial_s\Phi_u(s,x)du.$$
Hence
$$\partial_s\Phi_t(s,x)=-f(x)\exp\left(\int_s^tf'(\Phi_u(s,x))du\right).$$

 Denote  by $X=\{X_t, t\ge 0\}$ the Lévy process associated to the Poisson measure $N$ as in   (\ref{eq:lev-pr}) (obviously with $\gamma=\sigma=0$). The processes $X$ and $Z$
 jump at the same times, and the jumps height of the solution process
 $Z$ in a jump time $\tau$ is $h(\Delta X_\tau)$.
 Denote by $T_1<T_2<..$ the jump times of $N^{\Theta_m}$ (the dependence of $T_j$ on $m$ is suppressed to short the notations).
 The solution of (\ref{eq-approx}) in $\{N^{\Theta_m}_T=n\}$ is given by
\begin{align*}
 Z_t^{(m)}&=\Phi_t(0,x),\quad t\in[0,T_1),\\
 Z^{(m)}_{T_1}&=\Phi_{T_1}(0,x)+h(\Delta X_{T_1}),\\
 Z_t^{(m)}&=\Phi_{t}(T_1,Z^{(m)}_{T_1}),\quad t\in[T_1,T_2),
 \end{align*}
 and so on. See Figure \ref{grafic}.

 \begin{figure}[h]
\centering

\begin{tikzpicture}%[scale=2.5]%[rotate=20]%[scale=0.8]

\draw[line width=1, color=gris,->] (0,0)--(10,0); \draw[line
width=1, color=gris,->] (0,0)--(0,5); \draw[line width=2pt]
plot[smooth,tension=.8] coordinates { (0,1) (2,1.2) (3,.5) (4,2)};
\draw[line width=2pt] plot[smooth,tension=.8] coordinates {(4,4)
(5,4.5) (6,3.5)}; \draw[line width=2pt] plot[smooth,tension=.8]
coordinates {(6,2) (7,3) (8,2) (9,4)}; \draw(-0.2,1) node {$x$};
\draw[line width=0.5pt] (4,-0.1)--+(0,0.2); \draw[line width=0.5pt]
(6,-0.1)--+(0,0.2); \draw(4,-0.4) node {$T_{1}$}; \draw(6,-0.4) node
{$T_{2}$}; \draw(3,2.7) node {$h(\Delta X_{T_{1}})$}; \draw[line
width=1pt,dashed] (4,2)--(4,4); \filldraw(4,4) circle(0.1);
\filldraw[fill=white](4,2) circle(0.1); \filldraw(6,2) circle(0.1);
\filldraw[fill=white](6,3.5) circle(0.1); \draw(9.8,-0.4) node
{$t$};
\end{tikzpicture}
\caption{A trajectory of $Z_t^{(m)}$.} \label{grafic}
\end{figure}

By Corollary \ref{cor:derjt}, and induction over   $n$
\begin{eqnarray*}\lefteqn{
D_t^{\R_0}Z^{(m)}_{T_n}}\\
&=& f(Z^{(m)}_{T_n-}){k}({T_n},\Delta
X_{T_n})\left(\frac{T_n}{T}-\1_{[0,T_n]}(t)\right)\\
&&+\sum_{i=1}^{n-1}\exp\left(\int_{T_i}^{T_n}f'(Z^{(m)}_u)du\right)
\left(f(Z^{(m)}_{T_i-})- f(Z^{(m)}_{T_i})\right){k}({T_i},\Delta
X_{T_i})\left(\frac{T_i}{T}-\1_{[0,T_i]}(t)\right)
\end{eqnarray*}
on the set $\left\{N_T^{\Theta_m}=n\right\}$. Thus we have that
\begin{eqnarray}\label{eq:condomsol}\lefteqn{
D_t^{\R_0}\left(Z^{(m)}_T\1_{\{N_T^{\Theta_m}=n\}}\right)}\nonumber\\
&=&\1_{\{N_T^{\Theta_m}=n\}}\left[-f(Z^{(m)}_{T_n})\exp\left(\int_{T_n}^Tf'
(Z^{(m)}_u)du\right)k(T_n,\Delta
X_{T_n})\left(\frac{T_n}{T}-\1_{[0,T_n]}(t)\right)\right.
\nonumber\\
&&+\left. \exp\left(\int_{T_n}^T f'(Z^{(m)}_u)du\right)
D_t^{\R_0}Z^{(m)}_{T_n}\right]\nonumber\\
&=&\sum_{i=1}^{n}\exp\left(\int_{T_i}^Tf'(Z^{(m)}_u)du\right)
\left(f(Z^{(m)}_{T_i-}) -f(Z^{(m)}_{T_i})\right)\nonumber\\
&&\quad\quad\quad\times {k}({T_i},\Delta
X_{T_i})\left(\frac{T_i}{T}-\1_{[0,T_i]}(t)\right)
\1_{\{N_T^{\Theta_m}=n\}}\nonumber\\
&=&\1_{\{N_T^{\Theta_m}=n\}}\int_{(0,T]\times\Theta_m}
\exp\left(\int_{s}^Tf'(Z^{(m)}_u)du\right)
\nonumber\\
& & \qquad \times \left(f(Z^{(m)}_{s-})-
f(Z^{\Lambda_m}_{s})\right){k}(s,y)\left(\frac{s}{T}-\1_{[0,s]}(t)\right)
dN(s,y).
\end{eqnarray}

Now we can see that the random variable $Z_T$ is in the domain of
$D^{\R_0,{k}}$.
\begin{proposition}\label{prop:alfder}
Under the assumption of this section, we have that $Z_T\in {\rm
Dom}\, D^{\R_0}$.  Moreover,
\begin{equation}\label{eq:sigint}
 D_t^{\R_0}Z_T=\int_{[0,T]\times\R_0}
\exp\left(\int_{s}^Tf'(Z_u)du\right) \left(f(Z_{s-})-
f(Z_{s})\right){k}(s,y)\left(\frac{s}{T}-\1_{[0,s]}(t)\right)
dN(s,y).\end{equation}
\end{proposition}

\noindent\textit{Proof:} The result follows from Proposition
\ref{pro:convsol}, letting $m\rightarrow\infty$ and from the
dominated convergence theorem. That theorem can be applied due to
the integrand in the right-hand side of (\ref{eq:condomsol}) is
bounded by $M|h(y){k}(s,y)|$, for some $M>0$, since $f'$ is bounded,
and the jumps of the solution are related with $h$.  $\qquad
\square$

\subsubsection{Case with finite L\'evy measure}

In this section we analyze  the existence of a density for $Z_T$ in
the case that the L\'evy measure is finite and the drift is a
monotone function.
\begin{theorem}\label{pro:nou1}
Assume that $\nu(\R_0)<\infty$. Consider equation (\ref{eq:mondri}), where $f$ is a momotone function with
bounded derivative and $h\in L^2(\R_0,\nu))$ with
 $h(y)\neq0$ for $y\in\R_0$. Then $Z_T$ is absolutely continuous on
the set $\left\{N_T^{\R_0}\ge 1\right\}$.
\end{theorem}

\noindent\textit{Proof:} We first assume that $f$ is an increasing
function.  We choose for $k$ the function  $(-h\wedge 1)\vee (-1)$
that depends only on $y$. Then, $(f(Z_{s-})-f(Z_{s-}+h(\Delta
X_s)))k(s,\Delta X_s)$ is strictly positive in the jump points.
Thus, to finish the proof,  we only need to apply Proposition
\ref{pro:abcon} and observe that
$$\int_0^T\left(\frac{s}{T}-\1_{[0,s]}(t)\right)^2dt=s(1-\frac{s}{T}),\quad s\in[0,T],$$
and
$$\int_0^T\left(\frac{s}{T}-\1_{[0,s]}(t)\right)\left(\frac{r}{T}-\1_{[0,r]}(t)\right)dt=
s(1-\frac{r}{T}),\quad \textrm{for}\,  0\le s<r\le T.$$ Hence the
two integrals are positive for $s<t$.
 Indeed, we
have that
\begin{eqnarray*}\label{eq:qfalta}\lefteqn{
\int_0^T \left(D_t^{\R_0}Z_T\right)^2dt}\\
&=& \sum_{i=1}^{N_T^{\R_0}}
\exp\left(2\int_{T^{\R_0}_i}^Tf'(Z_u)du\right)
\left(f(Z_{T^{\R_0}_i}-) -f(Z_{T^{\R_0}_i})\right)^2 ({k}(\Delta
X_{T^{\R_0}_i}))^2T^{\R_0}_i\left(1-\frac{T^{\R_0}_i}{T}\right)\\
&&+2\sum_{1\le i<j\le
N_T^{\R_0}}\exp\left(\int_{T^{\R_0}_i}^Tf'(Z_u)du\right)
\left(f(Z_{T^{\R_0}_i}-) -f(Z_{T^{\R_0}_i})\right) {k}(\Delta
X_{T^{\R_0}_i})\\
&&\quad\quad\times \exp\left(\int_{T^{\R_0}_j}^Tf'(Z_u)du\right)
\left(f(Z_{T^{\R_0}_j}-) -f(Z_{T^{\R_0}_j})\right) {k}(\Delta
X_{T^{\R_0}_j})T^{\R_0}_i\left(1-\frac{T^{\R_0}_j}{T}\right),
\end{eqnarray*}
which is bigger that zero on the set $\left\{N_T^{\R_0}\ge
1\right\}$. Therefore $Z_T$ absolutely continuous in this set.

Finally we can proceed similarly in the case that $f$ is decreasing
using the function $k=(h\wedge 1)\vee (-1)$.$\quad\square$
\subsubsection{Case with infinite L\'evy measure}
Now we deal with the case that $f$  is only monotone on an
neighborhood of the initial condition $x$. This problem has been
analyzed in Nourdin and Simon \cite{Nu-Si-06} using a stratification
method.
\begin{theorem}Consider equation (\ref{eq:mondri}), where $f$ is a function with a bounded
derivative and $h\in L^2(\R_0,\nu)\cap L^1(\R_0,\nu)$. Assume that
$\nu(\R_0)=\infty$, $h(y)\neq0$ for $y\in\R_0$, and that $f$
 is  monotone  on a neighborhood of
the point $x$.  Then, the random variable $Z_T$  is absolutely
continuous.
\end{theorem}

 \noindent\textit{Proof:} Here it is convenient to write the function $k$ that we use in the derivative,
 $D^{\Lambda,k}$.  Proceeding as in Theorem \ref{pro:nou1}, we assume
 that there is $\varepsilon>0$ such that
 $f$ is increasing on $(x-\varepsilon, x+\varepsilon)$;
 the proof for $f$ decreasing is similar.

 Note that (\ref{eq:mondri}), together with the fact that $f$ has a bounded
 derivative, implies that there exists $M>0$ such that
$$|Z_t-x|\le
M\int_0^t|Z_s-x|ds+t|f(x)|+\int_{[0,t]\times\R_0}|h(y)|dN(s,y),\quad
t\in [0,T].$$ Consequently, Gronwall's  lemma leads to write
$$|Z_t-x|\le\left(t|f(x)|+\int_{[0,t]\times\R_0}|h(y)|dN(s,y)\right)e^{MT},
\quad t\in [0,T].$$
 Therefore, there is $t_0\in(0,T)$ such that
\begin{equation}\label{eq:cota}
|Z_t-x|\le\frac{\varepsilon}{2}+e^{MT}\int_{[0,t]\times\R_0}|h(y)|dN(s,y),\quad
t\in [0,t_0].\end{equation}

Now, for $t\in[0,t_0]$,  let
$$A_t=\left\{e^{MT}\int_0^t\int_{\R_0}|h(y)|dN(s,y)>\frac{\varepsilon}{2}\right\}$$
and   $k_t:[0,T]\rightarrow\R$ a function in $C^1((0,T))$ such that
$k_t(s)>0$ for $s\in[0,t)$, and $k_t(s)=0$ for $s\geq t$. Define
$k^{(t)}(s,y)=-k_t(s)((h(y)\wedge 1)\vee (-1))$. As the Lévy measure
is $\infty,$ there are, with probability 1,  jumps in $(0,t),$ and
from (\ref{eq:cota}) we deduce that  the jump sizes satisfy
$|Z_s-x|\le \varepsilon$ for $s\le t\le t_0$ on the set $A_t^c$.

Indeed, in this case we have that, on $A^c_t$,
$(f(Z_{s-})-f(Z_s))k^{(t)}(s,y)\ge 0,\quad \textrm{for}\  0\le s\le
T$, with $(f(Z_{s-})-f(Z_s))k^{(t)}(s,y)>0$ if $Z$ has a jump at
time $s\in(0,t)$. So, (\ref{eq:sigint}) gives
$$\int_0^T \left(D_u^{\R_0,k^{(t)}}Z_T\right)^2du>0\quad\textrm{in}\ \ A_t^c.$$
 As a consequence, by Proposition \ref{pro:abcon},
the condition $\lambda(B)=0$ implies that
$$P(\{Z_T\in B\}\cap A_t^c)=0.$$

Finally, we observe that $A_t^c\subset A_{t'}^c$ for $t'<t,$ and
using  the Markov inequality, we have

\begin{align*}
P(A_t)&\le
\frac{4e^{2MT}}{\varepsilon^2}E\left(\left(\int_{[0,t]\times\R_0}|h(y)|dN(s,y)\right)^2\right)\\
&\le\frac{8e^{2MT}}{\varepsilon^2}\,t\left(\int_{\R_0}|h(y)|^2\nu(dy)+
T\left(\int_{\R_0}|h(y)|\nu(dy)\right)^2\right).
\end{align*}
Hence we see that $\lim_{t\downarrow 0} P( A_t)=0. $

So $P(\cup_{t<t_0,\, t\in \Q}A_t^c)=1$, and we can conclude, as
$P(\{Z_T\in B\})=\lim_{t\downarrow 0} P(\{Z_T\in B\}\cap A_t^c),$
that $Z_T$ is absolutely continuous. $\quad\square$
\subsection{A pure discontinuous equation with non zero Wronskian}
In this section we assume that the L\'evy measure $\nu$ is finite
and  consider the stochastic differential equation
\begin{equation}\label{eq:CPeq}
Z_t=x+\int_0^tf(Z_s)ds+\int_0^t\int_{\R_0}h(y)g(Z_{s-})dN(s,y),\quad
t\in[0,T],
\end{equation}
where $x\in\R$, $h\neq 0$ is a bounded function in $L^2(\R_0,\nu)$,
and $f:\R\rightarrow\R$ and $g:\R\rightarrow\R$ are two bounded
functions with two bounded derivatives  and one bounded derivative,
respectively. The existence of a density for $Z_T$ was analyzed  by
Carlen and Pardoux \cite{Ca-Pa-90} in the case that the involved
L\'evy process is a Poisson process. In the remaining of this
section we also assume that
\begin{equation}\label{eq:wron}
|h(y)W(g,f)(x)|>\frac12
||f''||_{\infty}||h||_{\infty}^2||g||_{\infty}^2, \quad x\in\R \,
\, y\in\R_0, \end{equation}
 where $W(g,f)=g'f-f'g$ is the
Wronskian of $g$ and $f$.

As in Section \ref{sec:4.2}, we consider the flow $\{\Phi_t(s,x):
0\le s\le t\le T \ \textrm{and}\ x\in\R\}$ associated with equation
(\ref{eq:CPeq}) and the family of jump times $\{T_j: j\in\N\}$
of the L\'evy process $X$. Remember that now we are dealing with a
finite L\'evy measure in this part of the paper.

Note that the fact that
$Z_{T_n}=\Phi_{T_n}(T_{n-1},Z_{T_{n-1}})+h(\Delta
X_{T_{n}})g(Z_{T_{n-}})$ jointly  with $T_0=0$ and Corollary
\ref{cor:derjt} allow to utilize induction on $n$ to see that, for any
$k\in {\cal K}$,
\begin{eqnarray}\label{eq:prider}
D_t^{\R_0}Z_{T_n}&=&f(Z_{T_n-}){ k}(T_n,\Delta
X_{T_{n}})\left(\frac{T_n}{T}-\1_{[0,T_n]}(t)\right)\nonumber\\
&&+f(Z_{T_n-})h(\Delta X_{T_{n}})g'(Z_{T_n-}){k}(T_n,\Delta
X_{T_{n}})\left(\frac{T_n}{T}-\1_{[0,T_n]}(t)\right)\nonumber\\
&&+\sum_{i<n}\exp{\left(\int_{T_i}^{T_n}f'(Z_u)du\right)}r_t(T_i,\Delta
X_{T_{i}})\nonumber\\
&&+ \sum_{m=2}^n
\int_{S_m(T_m)}\exp{\left(\int_{s_1}^{T_n}f'(Z_u)du\right)}
h(y_m)g'(Z_{s_m-})\cdots
h(y_2)g'(Z_{s_2-})r_t(s_1,y_1)\nonumber\\
&&\quad\phantom{jjjjjjjjjjjjjjjjjjjjjjjjjjjjjjjj}\times
dN(s_m,y_m)\cdots dN(s_1,y_1),
\end{eqnarray}
where
$$r_t(s,y)={k}(s,y)\left(\frac{s}{T}-\1_{[0,s]}(t)\right)\left(f(Z_{s-})-f(Z_s)+f(Z_{s-})h(y)g'(Z_{s-})\right)$$
and $$S_m(r)=\{(s_1,y_1;\ldots;s_m,y_m)\in([0,T]\times\R_0)^m: 0\le
s_1<\cdots<s_m\le r\}.$$

Now we are ready to see that $Z_T$ is in the domain of $D^{\R_0}$.
\begin{proposition}
Let $Z$ be the solution of equation (\ref{eq:CPeq}). Then $Z_T\in
{\rm Dom}\, D^{\R_0}$.
\end{proposition}
\noindent\textit{Proof:} By Corollary \ref{cor:derjt} and
(\ref{eq:prider}), we have $Z_T\1_{[N^{\R_0}_T=n]}$ is in the domain
of $D^{\R_0}$ and
\begin{eqnarray*}
D_t^{\R_0}(Z_T\1_{[N_T^{\R_0}=n]})&=&D_t^{\R_0}
\left(\Phi_T(T_n,Z_{T_n})\1_{\{N_T^{\R_0}=n\}}\right)\\
&=&
\1_{[N_T^{\R_0}=n]}\left\{-f(Z_{T_n})\exp{\left(\int_{T_n}^Tf'(Z_u)du\right)}
{k}(T_n,\Delta
X_{T_n})\left(\frac{T_n}{T}-\1_{[0,T_n]}(t)\right)\right.\\
&&+\left.\exp{\left(\int_{T_n}^Tf'(Z_u)du\right)}D^{\R_0}(Z_{T_n})\right\}\\
&=&\1_{[N_T^{\R_0}=n]}\left\{\sum_{i\le
n}\exp{\left(\int_{T_i}^{T}f'(Z_u)du\right)}r_t(T_i,\Delta
X_{T_{i}})\right.\\
&&\left.+ \sum_{m=2}^n
\int_{S_m(T_n)}\exp{\left(\int_{s_1}^{T}f'(Z_u)du\right)}
k(y_m)g'(Z_{s_m-})\cdots
h(y_2)g'(Z_{s_2-})r_t(s_1,y_1)\right.\\
&&\left.\quad\phantom{\int_{S_m(T_n)}}\times dN(s_1,y_1)\cdots
dN(s_m,y_m)\right\}.
\end{eqnarray*}
Hence, the fact that there is a constant $C>0$ such that
$|r_t(s,y)|< C{k}(s,y)$ and
\begin{eqnarray*}\lefteqn{\left|
\int_{S_m(T_n)}\exp{\left(\int_{s_1}^{T}f'(Z_u)du\right)}\right.}\\
&&\quad\quad\left. \phantom{\int_{S_m(T_n)}}\times
h(y_m)g'(Z_{s_m-})\cdots h(y_2)g'(Z_{s_2-})r_t(s_1,y_1)
dN(s_1,y_1)\cdots dN(s_m,y_m)\right|\\
&&\leq\1_{\{m\leq n\}}\frac{C}{m!}\left(\int_0^T
\int_{\R_0}\left|{k}(s,y)\right|dN(s,y)\right)\left(\int_0^T
\int_{\R_0}\left|h(y)\right|dN(s,y)\right)^{m-1}
\end{eqnarray*}
implies that $Z_T$ is in the domain of $D^{\R_0}$ with
\begin{eqnarray}\label{eq:derCP}
D_t^{\R_0}Z_T&=&
\int_{[0,T]\times\R_0}\exp{\left(\int_{s}^{T}f'(Z_u)du\right)}r_t(s,y)dN(s,y)\nonumber\\
&&+ \sum_{m=2}^\infty
\int_{S_m(T)}\exp{\left(\int_{s_1}^{T}f'(Z_u)du\right)}
h(y_m)g'(Z_{s_m-})\cdots
h(y_2)g'(Z_{s_2-})r_t(s_1,y_1)\nonumber\\
&&\times dN(s_1,y_1)\cdots dN(s_m,y_m).
\end{eqnarray}
Thus the proof is complete. $\qquad \square$

Finally we can see that $Z_T$ has a density.
\begin{theorem}
Let $Z$ be the solution of equation (\ref{eq:CPeq}). Then $Z_T$ has
a density on the set $[N_T^{\R_0}>0]$.
\end{theorem}

 \noindent\textit{Proof:} For a rational number $p\in(0,T)$ and a positive integer $n\ge 1$,
 we introduce the set $$A_{p,n}=\{N_T=n\}\cap\{T_{n-1}\le p<T_n\},$$ with the convention $T_0=0.$
 We
 choose a function $h_p:[0,T]\rightarrow \R_+$ of class $C^1((0,T))$ such
 that $h_p(s)=0$ for $s\le p$, and $h_p(s)>0$ for $s>p$. Then, (\ref{eq:derCP})
 implies that in the set $A_{p,n}$ we only have to consider in the derivative  the jump at $T_n$,
\begin{eqnarray*}
 D_t^{\R_0,h_p}Z_T&=&\exp\left(\int^T_{T_n}f'(Z_u)du\right)h_p(T_n)\left(\frac{T_n}{T}-
 \1_{[0,T_n]}(t)\right)\\
 &&\quad\times
 \left(f(Z_{T_n-})-f(Z_{T_n})+f(Z_{T_n-})h(\Delta
 X_{T_n})g'(Z_{T_n-})\right),
\end{eqnarray*}
which is different than zero due to condition (\ref{eq:wron}). Hence
Proposition \ref{pro:abcon} implies that $P(\{Z_T\in B\}\cap
A_{p,n})=0$ for any Borel set $B\in{\cal B}(\R)$ such that
$\lambda(B)=0$. Thus, for this Borel set we have
$$P(\{Z_T\in B\}\cap\{N_T>0\})=P(\{Z_T\in B\}\cap(\cup_{n\in\N,p\in\Q}
A_{p,n}))=0,$$ and the proof is finished. $\qquad \square$
 \vskip .4cm
 \noindent \textbf{Acknowledgment}. The authors thank Cinvestav-IPN,
   Universitat
Aut\`onoma de Barcelona and Universitat de Barcelona   for their
hospitality and economical support during the realization of this
paper.

\end{document}